\documentclass[12pt]{article}
\usepackage{graphicx}
\usepackage{amsmath,amsthm,amssymb,enumerate}
\usepackage{euscript,mathrsfs}
\usepackage[left=2cm,right=2cm,top=3.5cm,bottom=3.5cm]{geometry}
\usepackage{color}
\catcode`\@=11 \@addtoreset{equation}{section}

\catcode`\@=12

\allowdisplaybreaks

\newtheorem{Theorem}{Theorem}[section]
\newtheorem{Proposition}[Theorem]{Proposition}
\newtheorem{Lemma}[Theorem]{Lemma}
\newtheorem{Corollary}[Theorem]{Corollary}

\theoremstyle{definition}
\newtheorem{Definition}[Theorem]{Definition}

\newtheorem{Remark}[Theorem]{Remark}

\newcommand{\bTheorem}[1]{
\begin{Theorem} \label{T#1} }
\newcommand{\eT}{\end{Theorem}}

\newcommand{\bProposition}[1]{
\begin{Proposition} \label{P#1}}
\newcommand{\eP}{\end{Proposition}}

\newcommand{\bLemma}[1]{
\begin{Lemma} \label{L#1} }
\newcommand{\eL}{\end{Lemma}}

\newcommand{\bCorollary}[1]{
\begin{Corollary} \label{C#1} }
\newcommand{\eC}{\end{Corollary}}

\newcommand{\bRemark}[1]{
\begin{Remark} \label{R#1} }
\newcommand{\eR}{\end{Remark}}

\newcommand{\bDefinition}[1]{
\begin{Definition} \label{D#1} }
\newcommand{\eD}{\end{Definition}}
\newcommand{\intON}[1]{\int_{\mathcal{T}^N} #1 \ \dx }

\newcommand{\Del}{\Delta_x}

\newcommand{\intOR}[1]{\int_{R^N} #1 \ \dx}

\newcommand{\bfphi}{\boldsymbol{\varphi}}

\newcommand{\bFormula}[1]{
\begin{equation} \label{#1}}
\newcommand{\eF}{\end{equation}}

\newcommand{\Ov}[1]{\overline{#1}}

\newcommand{\DC}{C^\infty_c}

\newcommand{\vr}{\varrho}
\newcommand{\vre}{\vr_\ep}

\newcommand{\vue}{\vu_\ep}

\newcommand{\vu}{\vc{u}}
\newcommand{\vc}[1]{{\bf #1}}

\newcommand{\Div}{{\rm div}_x}
\newcommand{\Grad}{\nabla_x}

\newcommand{\dx}{\,{\rm d} {x}}

\newcommand{\dt}{\,{\rm d} t }

\newcommand{\dxdt}{\dx \ \dt}
\newcommand{\intO}[1]{\int_{\Omega} #1 \ \dx}

\newcommand{\ep}{\varepsilon}

\definecolor{Cgrey}{rgb}{0.85,0.85,0.85}
\definecolor{Cblue}{rgb}{0.50,0.85,0.85}
\definecolor{Cred}{rgb}{1,0,0}
\definecolor{fancy}{rgb}{0.10,0.85,0.10}

\newcommand\Cbox[2]{%
    \newbox\contentbox%
    \newbox\bkgdbox%
    \setbox\contentbox\hbox to \hsize{%
        \vtop{
            \kern\columnsep
            \hbox to \hsize{%
                \kern\columnsep%
                \advance\hsize by -2\columnsep%
                \setlength{\textwidth}{\hsize}%
                \vbox{
                    \parskip=\baselineskip
                    \parindent=0bp
                    #2
                }%
                \kern\columnsep%
            }%
            \kern\columnsep%
        }%
    }%
    \setbox\bkgdbox\vbox{
        \color{#1}
        \hrule width  \wd\contentbox %
               height \ht\contentbox %
               depth  \dp\contentbox
        \color{black}
    }%
    \wd\bkgdbox=0bp%
    \vbox{\hbox to \hsize{\box\bkgdbox\box\contentbox}}%
    \vskip\baselineskip%
}

\date{}

\begin{document}


\title{On the low Mach number limit for the compressible Euler system}

\author{Eduard Feireisl
\thanks{The research of E.F.~leading to these results has received funding from the
European Research Council under the European Union's Seventh
Framework Programme (FP7/2007-2013)/ ERC Grant Agreement
320078. The Institute of Mathematics of the Academy of Sciences of
the Czech Republic is supported by RVO:67985840.} \thanks{This work was partially supported by the Simons - Foundation grant 346300 and the Polish Government MNiSW 2015-2019 matching fund.}
\and Christian Klingenberg\footnotemark[2]
\and Simon Markfelder\footnotemark[2] }

\date{May 24, 2017}

\maketitle

\bigskip

\centerline{Institute of Mathematics of the Academy of Sciences of the Czech Republic}
\centerline{\v Zitn\' a 25, CZ-115 67 Praha 1, Czech Republic}
\centerline{}
\centerline{Dept. of Mathematics, W\"urzburg University, Germany}

\medskip

\bigskip

\begin{abstract}
In this paper, we propose a new approach to singular limits of inviscid fluid flows based on the concept of 
dissipative measure--valued solutions. We show that dissipative measure-valued solutions of the compressible Euler equations converge to the smooth solution of the incompressible Euler system when the Mach number tends to zero. This holds both for well-prepared and ill-prepared initial data, where in the latter case the presence of acoustic waves causes difficulties. However this effect is eliminated on unbounded domains thanks to dispersion.

\end{abstract}


\tableofcontents

\section{Introduction}
\label{ii}

We propose a new approach to singular limits for inviscid fluid flows based on the concept of \emph{measure--valued solution} for the primitive system.
Specifically, we consider the barotropic compressible inviscid Euler equations in two and three space dimensions,
\begin{align} \label{compEuler1}
\partial_t \vr + \Div (\vr \vc{u}) &=0\\ \label{compEuler2}
\partial_t (\vr \vc{u}) + \Div (\vr \vc{u} \otimes \vc{u}) +  \Grad p(\varrho) &= 0
\end{align}
where $\varrho = \varrho(t,x)$ represents the mass density, $\vu = \vu(t,x)$  the velocity vector, and $p = p(\rho(x,t))$ the pressure. To avoid technicalities,
we focus on the iconic
example of the isentropic pressure--density state equation $p = a \varrho^\gamma$, with $\gamma > 1$, although more general cases can be treated as well.

One may rescale these equations by non-dimensionalization. After combining terms appropriately (setting the so--called Strouhal number equal to one) one reaches the following system
\begin{align} \label{compEulerEp1}
\partial_t \vre + \Div (\vre \vue) &=0\\ \label{compEulerEp2}
\partial_t (\vre \vue) + \Div (\vre \vue \otimes \vue) + \frac{1}{\ep^2} \Grad p(\vre) &= 0
\end{align}
where $\ep$ is called the Mach number. It represents the norm of the velocity divided by the sound speed.
For a more detailed derivation of this see the appendix in \cite{Barsukow2017} or Klainerman and Majda \cite{KM1}.
We consider the {asymptotic limit} of solutions $(\vre, \vue)$ for $\ep \rightarrow 0$. This process
represents a bridge between compressible and incompressible fluid flows. Indeed one can expand the dependent variables in terms of $\ep$. For example for the pressure we have
$$ p = p^{(0)} + \ep p^{(1)} + \ep^2 p^{(2)} + O(\ep^{3}) \quad .$$
We now collect terms of the same order and find that the zeroth and first order term in the pressure expansion are constant while
the zeroth order term of the velocity (which we shall call $\vc{v}$) satisfies the incompressibility condition $ \Div \vc{v} = 0$.
The resulting limiting equations are (setting the zeroth order term of density to be constant, and now calling the second order term pressure $p^{(2)} = \Pi $)
\begin{align}
\label{incompEuler1}
\Div \vc{v} &= 0 \\ \label{incompEuler2}
\partial_t \vc{v} + \vc{v} \cdot \Grad \vc{v} + \Grad \Pi &= 0 \quad .
\end{align}

The
initial data for the compressible equations for which the zeroth and first order term of the pressure are constant and the initial velocity
is solenoidal are called {\em well-prepared}. For the well-prepared data the above formal derivation has been made rigorous by e.g.  \cite{EB1}, \cite{KM1}, \cite{Uka}, \cite{Asano1987}, \cite{MESC1}. For a survey see \cite{SCHO2}. All these authors assume that the solutions of the compressible flow are smooth. However, as is well known, solutions of the compressible Euler system develop singularities in a finite time no matter how smooth and/or small the initial data are. One of the principal difficulties of this approach is therefore showing that the life span of the
classical solutions is in fact independent of the Mach number.

The hypothesis of smoothness of solutions is therefore quite restrictive and even not appropriate in the context of compressible inviscid fluids. On the other hand, the limit
incompressible Euler system, at least if considered in two space dimensions, admits global-in-time smooth solutions for smooth initial data. The existence of global smooth solutions for the incompressible Euler system in three space dimensions is an outstanding open problem.

To achieve global results, it is more convenient to consider the weak solutions of the compressible Euler system.
Recently, the theory
of convex integration produced a large number of global-in-time weak solutions basically for any regular initial data, however, ``most of them'' apparently violate
the basic energy inequality associated to the system, see e.g. Chiodaroli \cite{Chiod}, DeLellis and Sz\' ekelyhidi \cite{DelSze3}.
{In addition, there is also a non-void family of ``wild'' initial data} that give rise to infinitely many weak solutions satisfying many of the conventional
admissibility criteria, see Chiodaroli et al. \cite{ChiDelKre}, \cite{ChiKre}, DeLellis and Sz\' ekelyhidi \cite{DelSze3}. In spite of these results, the existence
of global--in--time \emph{admissible} weak solutions for \emph{arbitrary} {(possibly smooth)} initial data remains largely open for the compressible Euler system.

In this paper, we propose a new
approach  based on the concept of \emph{dissipative measure-valued (DMV) solution} recently developed by Gwiazda et al. \cite{GSWW}, \cite{FGSWW1}.
Roughly speaking, they are measure-valued solutions of the compressible Euler system satisfying an appropriate form of energy inequality, see Section
\ref{M}. The energy dissipation is expressed via a \emph{dissipation defect} that in turn dominates the concentration measures that may develop
in the field equations. The main advantage of this approach can be summarized as follows:
\begin{itemize}
\item The (DMV) solutions to the barotropic Euler system exist globally in time for any finite energy initial data. Indeed they can be identified as cluster points
of solutions to the Navier--Stokes system in the regime of vanishing viscosity, asymptotic limits of suitable numerical schemes as well as
limits of other suitable approximate problems, cf. Ne\v cas at al. \cite{MNRR}.

\item Although the (DMV) solutions are very general objects that are in general not uniquely determined by the initial data, the convergence is
unconditional as soon as the limit system admits a smooth solution.

\item Convergence holds for both well-prepared and ill-prepared initial data as long as the spatial domains allows dispersion of acoustic waves in the latter case, see Sections \ref{I} and \ref{Il}, respectively.

\end{itemize}
\noindent
Due to the low regularity of the DMV solutions, our method yields convergence in a very weak sense, specifically, in the sense of the strong topology on the space of probability measures.

The paper is organized as follows. After having introduced the necessary preliminary material in Section \ref{P}, we state our main result in Section \ref{PM}. Section \ref{I}
is devoted to the incompressible limit for well-prepared initial data under periodic boundary conditions. In Section \ref{Il} contains the proof of convergence for the ill-prepared data for the problem on the whole space $R^N$.

\section{Preliminaries and main result}
\label{P}

In this section, we collect some basic facts about (DMV) solutions and state our main result. The symbol $\Omega$ will denote the spatial domain occupied by the
fluid. We focus on two typical examples: Periodic boundary conditions, where $\Omega$ can be identified with
the ``flat'' torus
\[
\Omega = \mathcal{T}^N = \left( [-1,1] |_{\{ -1,1 \}} \right)^N,
\]
and $\Omega = R^N$, $N=2,3$.

\subsection{Measure--valued solutions to the compressible Euler system}
\label{M}

Let
\begin{equation} \label{M1}
\mathcal{Q} = \left\{ [\vr, \vc{m}] \ \Big| \ \vr \in [0, \infty), \ \vc{m} \in R^N \right\}
\end{equation}
be the natural phase space associated to solutions $[\vr, \vc{m}]=[\vr, \vr\vu]$ of the compressible Euler system \eqref{compEuler1}, \eqref{compEuler2}.

{A \emph{dissipative measure-valued (DMV) solution} to the compressible Euler system \eqref{compEuler1}, \eqref{compEuler2}
consists of a parameterized family of probability measures $Y_{t,x}$, $t \in (0,T)$, $x \in \Omega$,
\[
Y_{t,x} \in L^\infty_{{\rm weak}-(*)} \left((0,T) \times \Omega; \mathcal{P}(\mathcal{Q}) \right)
\]
and a non-negative function $\mathcal{D} \in L^\infty(0,T)$ called \emph{dissipation defect} satisfying:
\begin{itemize}
\item
{\bf Equation of continuity.}

\begin{equation} \label{M2}
\int_0^T \intO{ \left[ \left< Y_{t,x}; \vr \right> \partial_t \varphi + \left< Y_{t,x}; \vc{m} \right> \cdot \Grad \varphi \right] } \dt =
- \intO{\left< Y_{0,x}; \vr \right> \varphi(0, \cdot) } - \int_0^T \int_{{\Omega}} \Grad \varphi \cdot {\rm d} \mu^C_D
\end{equation}
for all $\varphi \in \DC([0,T) \times {\Omega})$ and a signed measure $\mu^C_D \in \mathcal{M} ([0,T] \times {\Omega}; R^N) $
called concentration defect.

\item
{\bf Momentum equation.}
\begin{equation} \label{M3}
\begin{split}
\int_0^T &\intO{ \left[ \left< Y_{t,x}; \vc{m} \right> \cdot \partial_t \bfphi + \left< Y_{t,x}; \frac{\vc{m} \otimes \vc{m}}{\vr} \right> : \Grad \bfphi
+ \left< Y_{t,x}; p(\vr) \right> \Div \bfphi
 \right] } \dt \\ &=
- \intO{\left< Y_{0,x}; \vc{m} \right> \cdot \bfphi(0, \cdot) } - \int_0^T \int_{\Ov{\Omega}} \Grad \bfphi : {\rm d} \mu^M_D
\end{split}
\end{equation}
for all $\bfphi \in \DC([0,T) \times {\Omega}; R^N)$ and a signed measure $\mu^M_D \in \mathcal{M} ([0,T] \times {\Omega}; R^{N \times N})$.

\item
{\bf Energy balance.}

The energy inequality
\begin{equation} \label{M5}
\begin{split}
\intO{ \left< Y_{\tau,x}; \frac{1}{2} \frac{ |\vc{m}|^2}{\vr} + P(\vr) - P'(\Ov{\vr})(\vr - \Ov{\vr}) - P(\Ov{\vr}) \right>} + \mathcal{D}(\tau) \\
\leq
\intO{\left< Y_{0,x}; \frac{1}{2} \frac{ |\vc{m}|^2}{\vr} + P(\vr) - P'(\Ov{\vr})(\vr - \Ov{\vr}) - P(\Ov{\vr})\right>}
\end{split}
\end{equation}
holds for a.a. $\tau \in (0,T)$ and a certain constant $\Ov{\vr} > 0$,  where $P(\vr) = \vr \int_1^\vr \frac{p(z)}{z^2} \ {\rm d}z$ called \emph{pressure potential}.
\item

{\bf Compatibility condition.}
\begin{equation} \label{M6}
\int_0^\tau \int_{{\Omega}} \left[ | \mu^C_D| + |\mu^M_D| \right] \dxdt \leq \int_0^\tau \xi (t) \mathcal{D}(t) \ \dt
\ \mbox{for a.a.}\ \tau \in [0,T], \ \xi \in L^1(0,T).
\end{equation}

\end{itemize}

\begin{Remark} \label{RM1}

Strictly speaking, the expressions containing the concentration defect in (\ref{M2}), (\ref{M3}) should be written
\[
\left< \mu^C_D; \Grad \varphi \right> \ \mbox{instead of}\
\int_0^T \int_{\Ov{\Omega}} \Grad \varphi \cdot {\rm d} \mu^C_D,\
\mbox{and}\ \left< \mu^C_M; \Grad \varphi \right> \ \mbox{instead of}\
\int_0^T \int_{\Ov{\Omega}} \Grad \varphi : {\rm d} \mu^M_D.
\]

Similarly, we should have written
\[
\left\| \mu^C_D \right\|_{[0,\tau] \times {\Omega} }
+ \left\| \mu^M_D \right\|_{[0,\tau] \times {\Omega} } \leq  \int_0^\tau \xi (t) \mathcal{D}(t) \ \dt,
\]
rather than
\[
\int_0^\tau \int_{{\Omega}} \left[ | \mu^C_D| + |\mu^M_D| \right] \dxdt \leq \int_0^\tau \xi (t) \mathcal{D}(t) \ \dt
\]
in (\ref{M6}).
\end{Remark}

\begin{Remark} \label{RM11}

In contrast with the original definition introduced in \cite{GSWW}, we prefer to work with the natural phase variable, namely the density $\vr$ and the
momentum $\vc{m} = \vr \vu$, similarly to \cite{FGSWW1}.

\end{Remark}

\begin{Remark} \label{RM12}

The constant $\Ov{\vr}$ in (\ref{M5}) can be taken arbitrary if $\Omega = \mathcal{T}^N$ and becomes relevant only for $\Omega = R^N$, where it represents
the far field limit of the density,
\[
\vr \to \Ov{\vr} \ \mbox{as}\ |x| \to \infty.
\]

\end{Remark}

\begin{Remark} \label{RM13}

The functions
\[
[ \vr, \vc{m}]  \mapsto \frac{\vc{m} \otimes \vc{m}}{\vr}, \ [\vr, \vc{m}] \mapsto \frac{ |\vc{m}|^2}{\vr}
\]
are singular at the boundary of the phase space $\mathcal{Q}$, namely on the vacuum zone $\vr = 0$. We set
\[
\frac{ |\vc{m}|^2 }{\vr} = \left\{ \begin{array}{l} 0 \ \mbox{if}\ \vr \geq 0, \ \vc{m} = 0\\ \infty
\ \mbox{if}\ \vr = 0, \ \vc{m} \ne 0    \end{array} \right.
\]
on the singular set. Accordingly, the function $[\vr, \vc{m}] \mapsto \frac{ |\vc{m}|^2 }{\vr}$ is convex lower semi-continuous on $\mathcal{Q}$.
Now it follows from the energy inequality (\ref{M5}) that $[\vr, \vc{m}] \mapsto \frac{ |\vc{m}|^2 }{\vr}$ is integrable with respect to
$Y_{t,x}$ for a.a. $t,x$. In particular,
\[
{\rm supp}[Y_{t,x}] \cap \left\{ [\vr, \vc{m} ]\in \mathcal{Q} \ \Big|\ \vr = 0, \ \vc{m} \ne 0 \right\} = \emptyset.
\]

\end{Remark}

\medskip

In applications, the parameterized family $Y_{t,x}$ is the Young measure generated by an oscillating sequence of approximate solutions
$[\vr, \vc{m} = \vr \vu]$, while the measure $Y_0$ is determined by the initial conditions. Note, however, there are measure--valued solutions to system \eqref{compEuler1}, \eqref{compEuler2} that \emph{are not}
generated by any sequence of weak solutions, see \cite{ChFKWi}.

The measures $\mu^C_D$, $\mu^M_D$ characterize the so-called concentration defect. There is a more precise characterization of these terms as soon as
a measure-valued solution is identified as a suitable limit of a family of weak solutions,
see Gwiazda et al. \cite{GSWW}. Then typically $\mu^C_D = 0$, while
$\mu^M_D$ is the Young measure associated to the so-called recession function corresponding to the quantity $\sqrt{\vr} u_i \sqrt{\vr} u_j + p(\vr) \delta_{i,j}$
in the sense of Alibert and Bouchitt\' e \cite{AliBou}. In such a case, the concentration defect $\mathcal{D}$ can be equally given in terms
of the recession function associated to the energy $\frac{1}{2} |\sqrt{\vr} \vu|^2 + P(\vr)$. These quantities satisfy the compatibility condition
(\ref{M6}) as soon as
\[
\limsup_{ \vr \to \infty} \frac{p(\vr)}{P(\vr)} \leq p_\infty < \infty,
\]
which implies
\begin{equation} \label{M6a}
p(\vr) \leq c(\Ov{\vr}, p_\infty) \left[ P(\vr) - P'(\Ov{\vr})(\vr - \Ov{\vr} ) - P(\Ov{\vr}) \right]
\end{equation}
for all $\vr$ large enough. Accordingly,
the function $\xi$ in (\ref{M6}) then can be taken constant depending \emph{only} on  $p_\infty$, $\Ov{\vr}$.

\begin{Remark} \label{RM2}

In the low Mach number limit problem studied below, the pressure takes the form $\frac{1}{\ep^2}p(\vr)$, while the associated pressure potential
reads $\frac{1}{\ep^2} P(\vr)$. In accordance with (\ref{M6a}), the measure-valued solutions
introduced by Gwiazda et al. \cite{GSWW}
will satisfy the compatibility condition
(\ref{M6}) \emph{uniformly} for $\ep \to 0$. The same remains true in the more general setting introduced in \cite{FGSWW1} and considered in the present
paper as long as the measure-valued solutions are generated by suitable family of functions, for which the concentration defect is characterized
as the difference between the weak-(*) limit in the sense of measures and the biting limit of nonlinear compositions, cf. \cite{FGSWW1}.

\end{Remark}

Finally, we remark that the \emph{existence} of the dissipative measure--valued solutions, at least for the iconic pressure law
$p(\vr) = a \vr^\gamma$, $\gamma \geq 1$ can be easily established by means of an artificial/physical viscosity approximation. Neustupa \cite{Neustup} constructed
a variant of the measure-valued solutions by considering a higher viscosity approximation to the Euler system in the spirit of the general theory of
multipolar fluids developed by Ne\v cas, \v Silhav\' y, and collaborators \cite{NESI}. In view of the nowadays available existence theory
for the barotropic Navier stokes system, the measure-valued solutions of the compressible Euler can easily be identified with the cluster points
for $\delta \to 0$
of a family of weak solutions $[\vr_\delta, \vc{m}_\delta]$ of the Navier--Stokes system:
\[
\partial_t \vr_\delta + \Div (\vr_\delta \vu_\delta) = 0,
\]
\[
\partial_t (\vr_\delta \vu_\delta) + \Div (\vr_\delta \vu_\delta \otimes \vu_\delta) + \Grad p_\delta(\vr_\delta) =
\delta \Delta \vu_\delta + \delta \Grad \Div \vu_\delta,
\]
\[
\begin{split}
&\intO{ \left[ \frac{1}{2} \vr_\delta |\vu_\delta|^2 + P_\delta (\vr_\delta) - P'_\delta(\Ov{\vr})(\vr_\delta - \Ov{\vr}) -
P_\delta (\Ov{\vr})  \right](\tau, \cdot) } +
\int_0^\tau \intO{ \delta \left[ |\Grad \vu_\delta|^2 + |\Div \vu_\delta |^2 \right] } \\ &\leq
\intO{ \left[ \frac{1}{2} \vr_0 |\vu_0|^2 + P_\delta(\vr_0) -P'_\delta(\Ov{\vr})(\vr_0 - \Ov{\vr}) -
P_\delta (\Ov{\vr}) \right] }, \ p_\delta(\vr) = p(\vr) + \delta \vr^\Gamma, \ \delta > 0.
\end{split}
\]
Indeed the existence of the weak solutions $[\vr_\delta, \vr_\delta \vu_\delta]$ is guaranteed by the theory of Lions \cite{LI4} for
$N=1, 2, 3$ at least if $\Gamma \geq \Gamma (N)$.
In view of Remark \ref{RM2}, the compatibility condition (\ref{M6}) {will be satisfied for a suitable constant $\xi$
independent of $\delta$}.

\subsection{Relative energy inequality}

For a parameterized family $Y_{t,x}$ of probability measure defined on the phase space (\ref{M1}),
we introduce \emph{the relative energy} functional
\begin{equation} \label{M7}
\begin{split}
\mathcal{E} &\left(\vr, \vc{m} \ \Big| \ r, \vc{U} \right) \\ &= \intO{
\left< Y_{t,x}; \frac{1}{2} \vr  \left| \frac{\vc{m}}{\vr} - \vc{U}(t,x) \right|^2 + P(\vr) - P'(r(t,x))(\vr - r(t,x)) - P(r(t,x)) \right>},
\end{split}
\end{equation}
where $\vc{U}$, $r$ are continuously differentiable ``test functions'', $\vc{U}$, $r - \Ov{\vr}$ compactly supported in $\Omega$, $r > 0$.

For all (DMV) solutions $\vr,\vc{m}$ of the compressible Euler system,
the following relation can be deduced from \eqref{M2}--\eqref{M5} , see \cite{GSWW}:
\begin{equation} \label{M8}
\begin{split}
&\left[ \mathcal{E} \left(\vr, \vc{m} \ \Big| \ r, \vc{U} \right) \right]_{t = 0}^{t = \tau} + \mathcal{D}(\tau) \\
&\leq
\int_0^\tau \intO{ \left[ \left< Y_{t,x}; \vr \vc{U}(t,x) - \vc{m} \right> \cdot \partial_t \vc{U} +
\left< Y_{t,x}; \vr \vc{U}(t,x) - \vc{m} ) \otimes \frac{\vc{m}}{\vr} \right> : \Grad \vc{U} - \left< Y_{t,x}; p(\vr) \right> \Div \vc{U}  \right] } \dt \\
&+ \int_0^\tau \intO{ \left[ \left< Y_{t,x} ; r(t,x) - \vr \right> \frac{1}{r} \partial_t p(r) -
\left< Y_{t,x} ; \vc{m} \right> \cdot \frac{1}{r} \Grad p(r)    \right] } \ \dt\\
&+ \int_0^\tau \int_{{\Omega}} \left( \frac{1}{2} \Grad |\vc{U}|^2 - \Grad P'(r) \right)\cdot {\rm d} \mu^C_D - \int_0^\tau \int_{{\Omega}} \Grad \vc{U}:{\rm d}\mu^M_D
\end{split}
\end{equation}
for any
\begin{equation} \label{M8a}
\vc{U}, r \in C^1([0,T] \times \Omega),\ r > 0, \ {\rm supp}[ \vc{U}], \ {\rm supp}[r - \Ov{\vr}]
\ \mbox{compact in}\ [0,T] \times \Omega.
\end{equation}

\begin{Remark} \label{RM15}

Note that compactness of the support of the test functions claimed in (\ref{M8a}) is irrelevant if $\Omega = \mathcal{T}^N$ -
a compact set.

\end{Remark}

\subsection{Solutions of the target system}

It is expected the low Mach number limit velocity $\vc{v}$ is described by the incompressible Euler system (\ref{incompEuler1}), (\ref{incompEuler2}).
Our approach leans essentially on the fact the limit field $\vc{v}$ is a smooth function. Referring to the classical result of Kato \cite{Kato72},
\cite{KaLai} we know that
(\ref{incompEuler1}), (\ref{incompEuler2}) admits a solution $\vc{v}$, unique in the class
\begin{equation} \label{regclass}
\vc{v} \in C([0, T_{\rm max}); W^{k,2}(\Omega;R^N)), \
\partial_t \vc{v},\ \partial_t \Pi,\ \Grad \Pi \in C([0, T_{\rm max}); W^{k-1,2}(\Omega;R^N)),
\end{equation}
for some $T_{\rm max} > 0$,
as soon as
\[
\vc{v}_0 \in W^{k,2}(\Omega;R^N)),\ k > \frac{N}{2} + 1, \ \Div \vc{v}_0 = 0.
\]
Moreover, the solution exists globally in time, meaning $T_{\rm max} = \infty$, if $N=2$.

\section{Main results}
\label{PM}

Let $\vr_{0,\ep} = \vr_\ep(0, \cdot)$, $\vu_{0, \ep} = \vu_\ep(0, \cdot)$ be the initial data for the rescaled system (\ref{compEulerEp1}), (\ref{compEulerEp2}).
We suppose that
\[
\frac{\vr_{0, \ep} - \Ov{\vr} }{\ep} \to s_0,\
\vc{u}_{0,\ep} \to {\vc{u}_0}
\]
in a certain sense specified in the forthcoming section.
We say that the initial data are
\begin{itemize}
\item {\emph{well-prepared}}
if $s_0 = 0$, $\vc{u}_0=\vc{v}_0$, $\Div \vc{v}_0 = 0$;
\item {\emph{ill-prepared}}
otherwise.
\end{itemize}

In the context of (DMV) solutions, where the the distribution of the initial data is determined by the measure $Y_{0,x}^\ep$, \emph{well-prepared initial
data} translates to
\begin{equation}\label{wp1}
\intO{
\left< Y_{0,x}^\ep; \frac{1}{2} \vr  \left| \frac{\vc{m}}{\vr} - \vc{v}_0(x) \right|^2 + \frac{1}{\ep^2} \Big( P(\vr) - P'(\Ov{\vr})(\vr - \Ov{\vr}) - P(\Ov{\vr}) \Big) \right>} \to 0 \ \mbox{as}\ \ep \to 0,
\end{equation}
for certain constant $\Ov{\vr} > 0$ and a solenoidal function $\vc{v}_0$.

If the initial data are given in terms of the functions $\vr_{0,\ep}$, $\vu_{0,\ep}$, meaning
\[
Y^{\ep}_{0,\ep} = \delta_{\vr_{0,\ep}(x), [\vr_{0,\ep}(x) \vu_{0,\ep}(x)]},
\]
(\ref{wp1}) follows as soon as
\[
\begin{split}
\frac{\vr_{0, \ep} - \Ov{\vr} }{\ep} \ \mbox{bounded in}\ L^\infty(\Omega),\ \frac{\vr_{0, \ep} - \Ov{\vr} }{\ep} \to 0 \ \mbox{in}\ L^1(\Omega),\
\vc{u}_{0,\ep} \to \vc{v}_0 \ \mbox{in} \ L^2(\Omega; R^N), \ \Div \vc{v}_0 = 0.
\end{split}
\]

Similarly, the initial data are \emph{ill-prepared} if
\begin{equation}\label{illp1}
\begin{split}
\intO{
\left< Y_{0,x}^\ep; \frac{1}{2} \vr  \left| \frac{\vc{m}}{\vr} - \vc{u}_0(x) \right|^2 + \frac{1}{\ep^2} \left( P(\vr) - P'\Big(\Ov{\vr} +
\ep s_0 \Big)\Big(\vr - \Ov{\vr} - \ep s_0\Big) - P\Big(\Ov{\vr} + \ep s_0\Big) \right) \right>} &\to 0 \\ \mbox{as}\ \ep &\to 0,
\end{split}
\end{equation}
for certain constant $\Ov{\vr} > 0$, $s_0 \in L^\infty \cap L^1(\Omega)$, and $\vc{u}_0 = \vc{v}_0 + \Grad \Phi_0$, $\Div \vc{v}_0  = 0$.
In terms of ``deterministic'' initial data $\vr_{0,\ep}$, $\vu_{0,\ep}$ this can be rephrased as
\[
\begin{split}
\frac{\vr_{0, \ep} - \Ov{\vr} }{\ep} \ \mbox{bounded in}\ L^\infty(\Omega),\ \frac{\vr_{0, \ep} - \Ov{\vr} }{\ep} \to s_0 \ \mbox{in}\ L^1(\Omega),\\
\vc{u}_{0,\ep} \to {\vc{u}_0 = \vc{v}_0 + \Grad \Phi_0} \ \mbox{in} \ L^2(\Omega; R^N), \ \Div {\vc{v}_0} = 0.
\end{split}
\]

\subsection{Main result for the well--prepared data}

We consider the rescaled compressible Euler system with the periodic boundary conditions, $\Omega = \mathcal{T}^N$, equipped with the well-prepared
initial data.
\begin{Theorem} \label{MRwp}
Let $p \in C^1(0, \infty) \cap C[0, \infty)$ satisfy $p'(\vr) > 0$ whenever $\vr > 0$.
Let $\Omega = \mathcal{T}^N$, $N=2,3$, and let $\left\{ Y^\ep_{t,x} \right\}_{t \in [0,T]; x \in \mathcal{T}^N}$, $\mathcal{D}^\ep$ be a family of (DMV) solutions of the rescaled
compressible Euler system (\ref{compEulerEp1}), (\ref{compEulerEp2}), satisfying the compatibility condition
(\ref{M6}) with $\xi$ independent of $\ep$. Let the initial data $Y^\ep_{0,x}$ be well--prepared, meaning (\ref{wp1}) holds for
$\Ov{\vr} > 0$ and $\vc{v}_0 \in W^{k,2}(\mathcal{T}^N; R^N)$, $k > \frac{N}{2} + 1$, $\Div \vc{v}_0 = 0$. Finally, suppose that $T < T_{\rm max}$, where
$T_{\rm max}$ denotes the life span of the solution to the incompressible Euler system (\ref{incompEuler1}), (\ref{incompEuler2}) endowed with the initial data $\vc{v}_0$.

Then
\[
\mathcal{D}^\ep \to 0 \ \mbox{in}\ L^\infty(0,T),
\]
\[
{\rm ess}\sup_{t \in (0,T)} \intON{
\left< Y_{t,x}^\ep; \frac{1}{2} \vr  \left| \frac{\vc{m}}{\vr} - \vc{v}(t,x) \right|^2 + \frac{1}{\ep^2} \Big( P(\vr) - P'(\Ov{\vr})(\vr - \Ov{\vr}) - P(\Ov{\vr}) \Big) \right>} \to 0
\]
as $\ep \to 0$,
where $\vc{v}$ is the solution of the incompressible Euler system (\ref{incompEuler1}), (\ref{incompEuler2}) with the initial data $\vc{v}_0$.

\end{Theorem}

Theorem \ref{MRwp} asserts that the probability measures $Y_{t,x}$ shrink to their expected value as $\ep \to 0$, where the latter are characterized
by the constant value $\Ov{\vr}$ for the density and the solution $\vc{v}$ of the incompressible system. The result is restricted to the life span of
$\vc{v}$ if $N=3$ and is global for $N = 2$. The required smoothness of $\vc{v}_0$ could possibly be slightly relaxed. The proof of Theorem
\ref{MRwp} is given in Section \ref{I} below.

\subsection{Main result for the ill--prepared data}

Convergence in the ill--prepared case is ``polluted'' by the presence of acoustic waves generated by the component $s_0$, $\Grad \Phi_0$ of the limit
data. To eliminate this effect, we consider the unbounded physical space $\Omega = R^N$, where dispersion annihilates acoustic phenomena at least on compact
sets.

To simplify presentation, we also assume that the concentration defect $\mu^C_D$ in the equation of continuity (\ref{M2}) vanishes. This is not a very severy restriction as it is always satisfied as long as the (DMV) solutions are obtained as a limit of a family of approximate solutions satisfying a suitable form of the energy balance.

\begin{Theorem} \label{MRillp}
Let $p \in C^1(0, \infty) \cap C[0, \infty)$ satisfy
\begin{equation} \label{pressure}
p'(\vr) > 0 \ \mbox{for all}\ \vr > 0, \ \limsup_{\vr \to \infty} \frac{p(\vr)}{P(\vr)} = P_\infty < \infty, \
\liminf_{\vr \to \infty} \frac{p(\vr)} {\vr^\gamma } \geq p_\infty > 0 \ \mbox{for some}\ \gamma > 1.
\end{equation}
Let $\Omega = R^N$, $N=2,3$, and let $\left\{ Y^\ep_{t,x} \right\}_{t \in [0,T]; x \in \mathcal{T}^N}$, $\mathcal{D}^\ep$ be a family of (DMV) solutions of the rescaled
compressible Euler system (\ref{compEulerEp1}), (\ref{compEulerEp2}), with $\mu^C_D = 0$, satisfying the compatibility condition
(\ref{M6}) with $\xi$ independent of $\ep$.  Let the initial data $Y^\ep_{0,x}$ be ill--prepared, meaning (\ref{illp1}) holds for
$\Ov{\vr} > 0$ and $s_0 \in W^{k,2} \cap W^{k,1}(R^N)$, $\vc{u}_0 = \vc{v}_0 \in W^{k,2} \cap W^{k,1} (R^N; R^N)$, $k > \frac{N}{2} + 2$.
Finally, suppose that $T < T_{\rm max}$, where
$T_{\rm max}$ denotes the life span of the solution to the incompressible Euler system (\ref{incompEuler1}), (\ref{incompEuler2}) endowed with the initial data $\vc{v}_0
= P[\vu_0]$, where $P$ denotes the standard Helmholtz projection onto the space of solenoidal functions.

Then
\[
\mathcal{D}^\ep \to 0 \ \mbox{in}\ L^\infty(0,T),
\]
\[
{\rm ess}\sup_{t \in (\delta,T)} \int_B
\left< Y_{t,x}^\ep; \frac{1}{2} \vr  \left| \frac{\vc{m}}{\vr} - \vc{v}(t,x) \right|^2 + \frac{1}{\ep^2} \Big( P(\vr) - P'(\Ov{\vr})(\vr - \Ov{\vr}) - P(\Ov{\vr}) \Big) \right> \ \dx \to 0
\]
as $\ep \to 0$, for any compact $B \subset R^N$ and any $0< \delta < T$,
where $\vc{v}$ is the solution of the incompressible Euler system (\ref{incompEuler1}), (\ref{incompEuler2}) with the initial data $\vc{v}_0$.
\end{Theorem}

Note that the required regularity of the data $s_0$, $\vu_0$ is higher than in Theorem \ref{MRwp}. Moreover, strong decay of $s_0$, $\vu_0$ is necessary as
$|x| \to \infty$. Convergence to the target system is only local, both in time and space. This is inevitable due to the presence of acoustic
waves. The proof of Theorem \ref{MRillp} will be done in Section \ref{Il}.

\section{Incompressible limit for well--prepared initial data}
\label{I}

In this section, we prove Theorem \ref{MRwp}. For $Y^{\ep}_{t,x}$ - the (DMV) solution of the rescaled system - we denote
\[
\mathcal{E}_\ep \left( \vr, \vc{m} \Big| \Ov{\vr}, \vc{v} \right)
= \intON{
\left< Y^\ep_{t,x}; \frac{1}{2} \vr  \left| \frac{\vc{m}}{\vr} - \vc{v}(t,x) \right|^2 + \frac{1}{\ep^2} \Big( P(\vr) - P'(\Ov{\vr})(\vr - \Ov{\vr}) - P(\Ov{\vr}) \Big) \right>}
\]
the relative entropy associated to $\Ov{\vr}$, $\vc{v}$.

\subsection{Relative energy inequality}

As the quantities $r = \Ov{\vr}$, $\vc{U} = \vc{v}$ enjoy the regularity required in (\ref{M8a}),
they can be used as test functions in the relative entropy inequality (\ref{M8}):
\begin{equation} \label{I7}
\begin{split}
\mathcal{E}_\ep &\left(\vr, \vc{m} \ \Big| \ \Ov{\vr}, \vc{v} \right)(\tau) + \mathcal{D}^\ep (\tau) \\ &\leq
\intON{
\left< Y^\ep_{0,x}; \frac{1}{2} \vr  \left| \frac{\vc{m}}{\vr} - \vc{v}_0(x) \right|^2 + \frac{1}{\ep^2} \Big( P(\vr) - P'(\Ov{\vr})(\vr - \Ov{\vr}) - P(\Ov{\vr}) \Big) \right>}  \\
&+
\int_0^\tau \intON{ \left[ \left< Y^\ep_{t,x}; \vr \vc{v}(t,x) - \vc{m} \right> \cdot \partial_t \vc{v} +
\left< Y^\ep_{t,x}; (\vr \vc{v}(t,x) - \vc{m} ) \otimes \frac{\vc{m}}{\vr} \right> : \Grad \vc{v}  \right] } \dt \\
&+ \int_0^\tau \int_{\mathcal{T}^N}  \frac{1}{2} \Grad |\vc{v}|^2 \cdot {\rm d} \mu^{C,\ep}_D - \int_0^\tau \int_{\mathcal{T}^N} \Grad \vc{v}:{\rm d}\mu^{M,\ep}_D.
\end{split}
\end{equation}

As the initial data are well-prepared, we get
\begin{equation} \label{I7a1}
\intON{
\left< Y^\ep_{0,x}; \frac{1}{2} \vr  \left| \frac{\vc{m}}{\vr} - \vc{v}_0(x) \right|^2 + \frac{1}{\ep^2} \Big( P(\vr) - P'(\Ov{\vr})(\vr - \Ov{\vr}) - P(\Ov{\vr}) \Big) \right>} \to 0 \ \mbox{as}\ \ep \to 0.
\end{equation}
In addition, since the compatibility condition (\ref{M6}) is satisfied uniformly with respect to $\ep$, we deduce
\begin{equation} \label{I7a2}
\int_0^\tau \int_{\mathcal{T}^N}  \frac{1}{2} \Grad |\vc{v}|^2 \cdot {\rm d} \mu^{C,\ep}_D - \int_0^\tau \int_{\mathcal{T}^N} \Grad \vc{v}:{\rm d}\mu^{M,\ep}_D
\leq c\left( \| \vc{v}_0 \|_{W^{k,2}} \right) \int_0^\tau \xi \mathcal{D}^\ep \ \dt
\end{equation}

In view of (\ref{I7a1}), (\ref{I7a2}), the conclusion of Theorem \ref{MRwp} follows by Gronwall's lemma as soon as we show
\begin{equation} \label{I7a3}
\begin{split}
\int_0^\tau & \intON{ \left[ \left< Y^\ep_{t,x}; \vr \vc{v}(t,x) - \vc{m} \right> \cdot \partial_t \vc{v} +
\left< Y^\ep_{t,x}; (\vr \vc{v}(t,x) - \vc{m} ) \otimes \frac{\vc{m}}{\vr} \right> : \Grad \vc{v}  \right] } \dt
\\ &\leq \omega(\ep) + c \int_0^\tau (1 + \xi) \left( \mathcal{E}_\ep \left(\vr, \vc{m} \ \Big| \ \Ov{\vr}, \vc{v} \right) + \mathcal{D}^\ep \right)
\dt,\ \omega(\ep) \to 0 \ \mbox{as}\ \ep \to 0.
\end{split}
\end{equation}

\subsection{Estimates}

Our goal is to show (\ref{I7a3}).

\subsubsection{Step 1 - convective term}

We start by writing
\[
\begin{split}
&\intON{ \left< Y^\ep_{t,x}; (\vr \vc{v}(t,x) - \vc{m} ) \otimes \frac{\vc{m}}{\vr} \right> : \Grad \vc{v} } \\&=
\intON{ \left< Y^\ep_{t,x}; (\vr \vc{v}(t,x) - \vc{m} ) \otimes \frac{\vc{m} - \vr \vc{v} }{\vr} \right> : \Grad \vc{v} }
+ \intON{ \left< Y^\ep_{t,x}; \vr \vc{v}(t,x) - \vc{m} \right> \cdot \vc{v} \cdot \Grad \vc{v} },
\end{split}
\]
where, obviously,
\[
\intON{ \left< Y^\ep_{t,x}; (\vr \vc{v}(t,x) - \vc{m} ) \otimes \frac{\vc{m} - \vr \vc{v} }{\vr} \right> : \Grad \vc{v} }
\leq c\left( \| \vc{v}_0 \|_{W^{k,2}} \right) \mathcal{E}_\ep \left(\vr, \vc{m} \ \Big| \ \Ov{\vr}, \vc{v} \right).
\]
Moreover, as $\vc{v}$ fulfills equation (\ref{incompEuler2}), we may go back to (\ref{I7a3}) to deduce that (\ref{I7a3}) reduces to showing
\begin{equation} \label{I7a4}
\int_0^\tau  \intON{ \left< Y^\ep_{t,x}; \vc{m} - \vr \vc{v}(t,x) \right> \cdot \Grad \Pi } \dt
\leq \omega(\ep) + c \int_0^\tau (1 + \xi) \left( \mathcal{E}_\ep \left(\vr, \vc{m} \ \Big| \ \Ov{\vr}, \vc{v} \right) + \mathcal{D}^\ep \right)
\dt.
\end{equation}

\subsubsection{Step 2 - pressure estimates}

To see (\ref{I7a4}), we
deduce from {\eqref{M2}} that
\begin{equation} \label{I9a}
\begin{split}
\int_0^\tau &\intON{ \left< Y^\ep_{t,x}; \vc{m} \right> \cdot \Grad \Pi }\dt \\ &= - \int_0^\tau \intON{\left< Y^\ep_{t,x}; \vr \right>  \partial_t \Pi }\dt
+ \left[ \intON{ \left< Y^\ep_{t,x}; \vr \right> \Pi } \right]_{t = 0}^{t = \tau} { - \int_0^\tau \int_{\mathcal{T}^N} \Grad \Pi \cdot {\rm d} \mu^{C,\ep}_D } \\&=
- \int_0^\tau \intON{\left< Y^\ep_{t,x}; \vr - \Ov{\vr} \right>  \partial_t \Pi }\dt
+ \left[ \intON{ \left< Y^\ep_{t,x}; \vr - \Ov{\vr} \right> \Pi } \right]_{t = 0}^{t = \tau} { - \int_0^\tau \int_{\mathcal{T}^N} \Grad \Pi \cdot {\rm d} \mu^{C,\ep}_D }\\&=
- \ep \int_0^\tau \intON{\left< Y^\ep_{t,x}; \frac{\vr - \Ov{\vr}}{\ep} \right>  \partial_t \Pi }\dt
+ \ep \left[ \intON{ \left< Y^\ep_{t,x}; \frac{\vr - \Ov{\vr}}{\ep} \right> \Pi } \right]_{t = 0}^{t = \tau} \\  &- \int_0^\tau \int_{\mathcal{T}^N} \Grad \Pi \cdot {\rm d} \mu^{C,\ep}_D .
\end{split}
\end{equation}

Similarly, we may use the incompressibility condition $\Div \vc{v} = 0$ to obtain
\begin{equation} \label{I9b}
\int_0^\tau \intON{ \left< Y^\ep_{t,x}; \vr \vc{v}(t,x) \right> \cdot \Grad \Pi } \dt
= \ep \int_0^\tau \intON{ \left< Y^\ep_{t,x}; \frac{\vr - \Ov{\vr}}{\ep}  \right> \vc{v} \cdot \Grad \Pi } \dt
\end{equation}

Now observe that the rightmost integral in (\ref{I9a}) can be controlled by the dissipation defect $\mathcal{D}^\ep$. Consequently,
as the pressure $\Pi$ belongs to the regularity class (\ref{regclass}), in particular $\Pi$, $\partial_t \Pi$ and $\Grad \Pi$ are bounded continuous in
$[0,T] \times \mathcal{T}^N$, it is enough to establish a uniform bound
\begin{equation} \label{I10a}
\intON{\left< Y^\ep_{t,x}; \left| \frac{\vr - \Ov{\vr}}{\ep} \right| \right> } \leq c.
\end{equation}

\subsubsection{Step 3 - energy estimates}

As the (DMV) solutions satisfy the energy inequality (\ref{M5}), we deduce from (\ref{wp1}) that
\begin{equation} \label{I10b}
\frac{1}{\ep^2} \intON{ \left< Y^\ep_{t,x}; P(\vr) - P'(\Ov{\vr})(\vr - \Ov{\vr}) - P(\Ov{\vr}) \right> } \leq c
\ \mbox{uniformly as}\ \ep \to 0.
\end{equation}
Since
\[
P''(\vr) = \frac{p'(\vr)}{\vr} \ \mbox{for}\ \vr > 0,
\]
the function $P$ is strictly convex, and, consequently
\begin{equation} \label{I10c}
| \vr - \Ov{\vr} |^2 \leq c(\delta) \Big( P(\vr) - P'(\Ov{\vr})(\vr - \Ov{\vr}) - P(\Ov{\vr}) \Big)
\ \mbox{whenever}\ 0 < \delta \leq \vr, \Ov{\vr} \leq \frac{1}{\delta},\ \delta > 0,
\end{equation}
and
\begin{equation} \label{I10d}
\begin{split}
1 + | \vr - \Ov{\vr} | &+ P(\vr)  \leq c(\delta) \Big( P(\vr) - P'(\Ov{\vr})(\vr - \Ov{\vr}) - P(\Ov{\vr}) \Big)
\\ &\mbox{if}\ 0 < 2 \delta < \Ov{\vr} < \frac{1}{2 \delta}, \ \vr \in [0, \delta) \cup [\frac{1}{\delta}, \infty),\
\delta > 0.
\end{split}
\end{equation}

Combining (\ref{I10c}), (\ref{I10d}) with (\ref{I10b}) we obtain (\ref{I10a}). Theorem \ref{MRwp} has been proved.

\section{Incompressible limit for ill--prepared initial data}
\label{Il}

Our goal is to prove Theorem \ref{MRillp}. To begin, we introduce a function
$\chi = \chi(\vr)$ such that
\[
\chi(\vr) \in \DC(0, \infty),\ 0 \leq \chi \leq 1, \chi(\vr) = 1 \ \mbox{if}\ \frac{\Ov{\vr}}{2} \leq \vr \leq 2 \Ov{\vr}.
\]
For a function $H = H(\vr, \vc{m})$ we set
\[
H_{\rm ess}(\vr, \vc{m}) = \chi(\vr) H(\vr, \vc{m}),\
H_{\rm res}(\vr, \vc{m}) = (1 - \chi(\vr)) H(\vr, \vc{m}).
\]

\subsection{Energy bounds}

As the initial distribution $Y^\ep_{0,x}$ is ill--prepared, meaning satisfies (\ref{illp1}), and the the functions $s_0$, $\vu_0$ belong to
$L^\infty \cap L^1(R^N)$, the initial energy
\[
\intOR{ \left< Y^{\ep}_{0,x}; \frac{1}{2} \frac{ |\vc{m}|^2}{\vr} + \frac{1}{\ep^2} \Big( P(\vr) - P'(\Ov{\vr}) (\vr - \Ov{\vr}) - P(\Ov{\vr}) \Big) \right>}
\leq E_0
\]
is bounded uniformly for $\ep \to 0$. In accordance with the energy inequality we obtain
\begin{equation} \label{E1}
{\rm ess} \sup_{t \in (0,T)}
\intOR{ \left< Y^{\ep}_{t,x}; \frac{1}{2} \frac{ |\vc{m}|^2}{\vr} + \frac{1}{\ep^2} \Big( P(\vr) - P'(\Ov{\vr}) (\vr - \Ov{\vr}) - P(\Ov{\vr}) \Big) \right>}
\leq E_0.
\end{equation}
Thus, using estimates (\ref{I10c}), (\ref{I10d}), we may infer that
\begin{equation} \label{E2}
{\rm ess} \sup_{t \in (0,T)} \intOR{ \left< Y^{\ep}_{t,x}; \left| \left[ \frac{  \vr - \Ov{\vr} }{\ep} \right]_{\rm ess} \right|^2 \right> }
+ {\rm ess} \sup_{t \in (0,T)} \intOR{ \left< Y^{\ep}_{t,x}; \left[ \frac{P(\vr) + 1}{\ep^2} \right]_{\rm res}  \right> } \leq c
\end{equation}

Furthermore, we get
\begin{equation} \label{Il3b}
\intOR{ \left< Y^\ep_{\tau,x}; \left[ |\vc{m}|^2\right]_{\rm ess} \right> } \leq c
\intOR{ \left< Y^\ep_{\tau,x}; \left[ \vr \left|\frac{\vc{m}}{\vr}\right|^2 \right]_{\rm ess} \right> } \leq E_0.
\end{equation}
Seeing that
\[
[ |\vc{m}| ]_{\rm res} \leq \frac{|\vc{m}|}{\sqrt{\vr}}  \left[ \sqrt{\vr} \right]_{\rm res}
\]
we deduce
\[
\left[ |\vc{m}|^{\frac{2 \gamma}{\gamma + 1} } \right]_{\rm res} \leq c \left( \ep \vr \frac{ |\vc{m} |^2 }{\vr} +
\frac{1}{\ep} \left[ \vr^{\gamma} \right]_{\rm res} \right);
\]
whence
\begin{equation} \label{Il3e}
\intOR{ \left< Y^\ep_{\tau, x} ; \left[ |\vc{m}|^{\frac{2 \gamma}{\gamma + 1} } \right]_{\rm res} \right> } \leq \ep E_0.
\end{equation}

Finally, we recall Jensen's inequality
\begin{equation} \label{Jensen}
\left< Y^\ep_{\tau,x}; |\vc{F}| \right>^q \leq \left< Y^\ep_{\tau,x} |\vc{F} |^q \right>,\ q \geq 1.
\end{equation}

Consequently, the estimates \eqref{Il3b}--\eqref{Il3e} give rise to
\begin{equation} \label{Il3f}
\begin{split}
\left< Y^\ep_{\tau, \cdot}; \vc{m} \right> \ &\mbox{bounded in}\ \left[ L^2 + L^{\frac{2 \gamma}{\gamma + 1}} \right] (R^N, R^N),\\
\left< Y^\ep_{\tau, \cdot}; \left[ \frac{\vr - \Ov{\vr}}{\ep} \right]_{\rm ess} \right>\ &\mbox{bounded in}\ L^2 (R^N),
\\ {\ep^{- \frac{2}{\gamma}} } \left< Y^\ep_{\tau, \cdot}; \left[ \vr \right]_{\rm res} \right>\  &\mbox{bounded in}\ L^\gamma (R^N).
\end{split}
\end{equation}

\subsection{Acoustic equation}

Write
\[
\vc{u}_0 = \vc{v}_0 + \Grad \Phi_0, \ \vc{v}_0 = P[\vc{u}_0].
\]
The evolution of acoustic waves is described by the \emph{acoustic equation}
\begin{align}
\label{A1}
\ep \partial_t s_\ep + \Div (\Ov{\vr} \Grad \Phi_\ep ) &= 0 \\
\label{A2}
\ep \partial_t \Grad \Phi_\ep + \frac{p'(\Ov{\vr})}{\Ov{\vr}} \Grad s_\ep &= 0 \\ \nonumber
s(0, \cdot) = s_0, \ \Grad \Phi_\ep (0, \cdot) &= \Grad \Phi_0
\end{align}
considered in the whole space $R^N$, $N=2,3$.

\subsubsection{Acoustic energy}

Solutions of (\ref{A1}), (\ref{A2}) conserve the total (acoustic) energy, specifically,
\begin{equation} \label{Il2bis}
\frac{{\rm d}}{{\rm d}t} \intOR{ \left[ p'(\Ov{\vr}) s^2_\ep + \Ov{\vr}^2 |\Grad \Phi_\ep |^2 \right] } = 0.
\end{equation}
Differentiating the (linear) system (\ref{A1}), (\ref{A2}) we easily extend (\ref{Il2bis}) to
\begin{equation} \label{Il2}
\left\| s_\ep (\tau, \cdot) \right\|^2_{W^{k,2}(R^N)} + \left\| \Grad \Phi_\ep (\tau, \cdot) \right\|^2_{{W^{k,2}(R^N;R^N)}}
\leq c \left[ \left\| s_0 \right\|^2_{W^{k,2}(R^N)} + \left\| \Grad \Phi_0  \right\|^2_{{W^{k,2}(R^N;R^N)}} \right]
\end{equation}
for any $\tau \geq 0$, $k \geq 0$.

\subsubsection{Dispersion}

Finally, we report the dispersive estimates
\begin{equation} \label{Il3}
\begin{split}
\left\| s_\ep (\tau, \cdot) \right\|^2_{L^p(R^3)} &+ \left\| \Grad \Phi_\ep (\tau, \cdot) \right\|^2_{L^p (R^3;R^3)}\\
\\ &\leq c \left( 1 + \frac{\tau}{\ep} \right)^{(N-1) \left( \frac{1}{p} - \frac{1}{q} \right)}\left[ \left\| s_0 \right\|^2_{W^{k,q}(R^3)} + \left\| \Grad \Phi_0  \right\|^2_{{W^{k,q}(R^3;R^3)}} \right],
\end{split}
\end{equation}
$k \geq N \left( \frac{1}{q} - \frac{1}{p} \right)$, $2 \leq p \leq \infty$, $\frac{1}{p} + \frac{1}{q} = 1$, see Strichartz \cite{Strich}.

\subsection{Relative energy inequality}

The first observation is that
\[
r = \Ov{\vr} + \ep s_\ep , \vc{U} = \vc{v} + \Grad \Phi_\ep,
\]
where $\vc{v}$ is the solution of the incompressible Euler system (\ref{incompEuler1}), (\ref{incompEuler2}), and
$s_\ep$, $\Phi_\ep$ solve the acoustic system (\ref{A1}), (\ref{A2}) can be taken as test functions in the relative energy inequality
(\ref{M8}). Note that, strictly speaking, these functions do not belong to the class (\ref{M8a}) but decay sufficiently fast to their far field limit. Validity of (\ref{M8}) can be verified by a density argument. Note the the most problematic term containing the pressure can be handled as follows:
\[
\intOR{ \left< Y^\ep_{t,x}; p(\vr) \right> \Div \vc{U} } = \intOR{ \left< Y^\ep_{t,x}; p(\vr) - p(\Ov{\vr}) \right> \Div \vc{U} }
= \intOR{ \left< Y^\ep_{t,x}; p(\vr) - p(\Ov{\vr}) \right> \Del \Phi }.
\]

Writing
\[
\begin{split}
\mathcal{E}_\ep &\left(\vr, \vc{m} \ \Big| \ \Ov{\vr} + \ep s_\ep , \vc{v} + \Grad \Phi_\ep \right)
\\ &= \intOR{
\left< Y^\ep_{t,x}; \frac{1}{2} \vr  \left| \frac{\vc{m}}{\vr} - \vc{v}(t,x) - \Grad \Phi_\ep(t,x) \right|^2 \right> } \\  &+
\frac{1}{\ep^2} \intOR{ \left< Y^\ep_{t,x};  P(\vr) - P'(\Ov{\vr} +
\ep s_\ep(t,x) )(\vr - \Ov{\vr} - \ep s_{\ep}(t,x)) - P\left(\Ov{\vr} + \ep s_\ep(t,x)  \right) \right>},
\end{split}
\]
we obtain the relative energy inequality in the form
\begin{equation} \label{Il4}
\begin{split}
&\left[ \mathcal{E}_\ep \left(\vr, \vc{m} \ \Big| \ \Ov{\vr} + \ep s_\ep , \vc{v} + \Grad \Phi_\ep \right) \right]_{t = 0}^{t = \tau} + \mathcal{D}^\ep(\tau)
\leq \omega(\ep) \\
&+
\int_0^\tau \intOR{ \left[ \left< Y^\ep_{t,x}; \vr \vc{U}(t,x) - \vc{m} \right> \cdot \partial_t \vc{U} +
\left< Y^\ep_{t,x}; (\vr \vc{U}(t,x) - \vc{m} ) \otimes \frac{\vc{m}}{\vr} \right> : \Grad \vc{U} \right] }\dt
\\
&
- \frac{1}{\ep^2} \int_0^\tau \intOR{ \left< Y^\ep_{t,x}; p(\vr) - p(\Ov{\vr}) \right> \Delta \Phi_\ep   } \dt \\
&+ \frac{1}{\ep} \int_0^\tau \intOR{ \left[ \left< Y^\ep_{t,x} ; r(t,x) - \vr \right> P''(r) \partial_t s_\ep -
\left< Y^\ep_{t,x} ; \vc{m} \right> \cdot P''(r) \Grad s_\ep    \right] } \ \dt\\
& - \int_0^\tau \int_{R^N} \Grad \vc{U}:{\rm d}\mu^{M,\ep}_D,\ \omega(\ep) \to 0 \ \mbox{as}\ \ep \to 0.
\end{split}
\end{equation}

In view of the dispersive estimates (\ref{Il3}), the conclusion of Theorem \ref{MRillp} follows as soon as we show that the expression on the right--hand
side of (\ref{Il4}) vanishes for $\ep \to 0$. Similarly to the previous section, we use a Gronwall type arguments proceeding in several steps.

\subsubsection{Step 1 - convective term I}

Similarly to Section \ref{I}, one may use the compatibility condition (\ref{M6}) to control the error term
\[
\int_0^\tau \int_{R^N} \Grad \vc{U}:{\rm d}\mu^{M,\ep}_D \leq \| \Grad \vc{U} \|_{L^\infty} \int_0^\tau \xi (t) \mathcal{D}^\ep (t) \ \dt.
\]

Next, exactly
as in the well--prepared case, we write
\[
\begin{split}
&\left< Y^\ep_{t,x}; (\vr \vc{U}(t,x) - \vc{m} ) \otimes \frac{\vc{m}}{\vr} \right> : \Grad \vc{U} \\
&= \left< Y^\ep_{t,x}; (\vr \vc{U}(t,x) - \vc{m} ) \otimes \left( \frac{\vc{m}}{\vr} - \vc{U} \right) \right> : \Grad \vc{U}
+ \left< Y^\ep_{t,x}; \vr \vc{U}(t,x) - \vc{m} \right> \cdot \vc{U} \cdot \Grad \vc{U},
\end{split}
\]
to deduce that (\ref{Il4}) reduces to
\begin{equation} \label{Il5}
\begin{split}
&\mathcal{E}_\ep \left(\vr, \vc{m} \ \Big| \ \Ov{\vr} + \ep s_\ep , \vc{v} + \Grad \Phi_\ep \right)  (\tau)  + \mathcal{D}^\ep (\tau)
\leq \omega(\ep)
\\
&+
\int_0^\tau \intOR{ \left[ \left< Y^\ep_{t,x}; \vr \vc{U} - \vc{m} \right> \cdot \left( \partial_t \vc{U} + \vc{U} \cdot \Grad \vc{U} \right)
\right] }\dt
\\
&
- \frac{1}{\ep^2} \int_0^\tau \intOR{ {\left< Y^\ep_{t,x}; p(\vr) - p(\Ov{\vr}) \right> }\Delta \Phi_\ep   } \dt \\
&+ \frac{1}{\ep} \int_0^\tau \intOR{ \left[ \left< Y^\ep_{t,x} ; r - \vr \right> P''(r) \partial_t s_\ep {-
\left< Y^\ep_{t,x} ; \vc{m}  \right> } \cdot P''(r) \Grad s_\ep    \right] } \ \dt\\
& + c \int_0^\tau (1 + \xi(t)) \left[ \mathcal{E}_\ep \left(\vr, \vc{m} \ \Big| \ \Ov{\vr} + \ep s_\ep , \vc{v} + \Grad \Phi_\ep \right) + \mathcal{D}^\ep
\right] \dt
\end{split}
\end{equation}

\subsubsection{Step 2 - convective term II}

Next, we rewrite
\[
\begin{split}
\int_0^\tau &\intOR{ \left[ \left< Y^\ep_{t,x}; \vr \vc{U} - \vc{m} \right> \cdot \left( \partial_t \vc{U} + \vc{U} \cdot \Grad \vc{U} \right)
\right] }\dt\\ &= \int_0^\tau \intOR{ \left[ \left< Y^\ep_{t,x}; \vr \vc{U} - \vc{m} \right> \cdot \left( \partial_t \vc{v} + \vc{v} \cdot \Grad \vc{v} \right)
\right] }\dt \\ &+ \int_0^\tau \intOR{ \left[ \left< Y^\ep_{t,x}; \vr \vc{U} - \vc{m} \right> \cdot \left( \partial_t \Grad \Phi_\ep \right)
\right] }\dt \\ &+ \int_0^\tau \intOR{ \left[ \left< Y^\ep_{t,x}; \vr \vc{U} - \vc{m} \right> \cdot \Grad \Phi_\ep {\cdot} \Grad \vc{v}
\right] }\dt \\&+
\int_0^\tau \intOR{ \left[ \left< Y^\ep_{t,x}; \vr \vc{U} - \vc{m} \right> \otimes \vc{v} \right] : \nabla^2_x \Phi_\ep  }\dt
\\ &+ \frac{1}{2} \int_0^\tau \intOR{ \left[ \left< Y^\ep_{t,x}; \vr \vc{U} - \vc{m} \right> {\cdot} \Grad \left| \Grad \Phi_\ep \right|^2  \right]}\dt.
\end{split}
\]

First observe that the last three integrals can be controlled in terms of
\[
\int_0^\tau \left\| \Grad \Phi_\ep \right\|_{W^{1,p}(R^3; R^3)}\ {\dt}
\ \ \ \mbox{for some}\ p > 2 \ \mbox{sufficiently large},
\]
and, consequently, in accordance with the dispersive estimates \eqref{Il3} vanish in the asymptotic limit $\ep \to 0$. Indeed the desired estimates on
$\left< Y; \vr \right>$, $\left< Y; \vc{m} \right>$ follow form \eqref{Il3f}, while $\vc{v}$ is bounded being a smooth solution of the incompressible Euler system.

Next, we have
\[
\begin{split}
\int_0^\tau &\intOR{ \left[ \left< Y^\ep_{t,x}; \vr \vc{U} - \vc{m} \right> \cdot \left( \partial_t \vc{v} + \vc{v} \cdot \Grad \vc{v} \right)
\right] }\dt \\&=
\int_0^\tau \intOR{ \left< Y^\ep_{t,x}; \vc{m} \right> \cdot \Grad \Pi }\dt - \int_0^\tau \intOR{ \left< Y^\ep_{t,x}; \vr \right> \vc{U} \cdot \Grad \Pi } \dt,
\end{split}
\]
where the former term on the right--hand side may be handled exactly as in (\ref{I9a}). As for the latter, we get
\[
\begin{split}
&\left| \int_0^\tau \intOR{ \left< Y^\ep_{t,x}; \vr \right>  \vc{U} \cdot \Grad \Pi } \dt \right| \\ \leq &\ep \left|
\int_0^\tau \intOR{ \left< Y^\ep_{t,x}; \frac{\vr - \Ov{\vr}}{\ep} \right>  \vc{U} \cdot \Grad \Pi } \dt \right| + \Ov{\vr} \left|
\int_0^\tau \intOR{ \Grad \Phi_\ep \cdot \Grad \Pi } \dt \right|,
\end{split}
\]
where the first term is small because of \eqref{Il3f}, while the second one vanishes for $\ep \to 0$ because of dispersive estimates. Indeed the pressure
$\Pi$ may be computed by means of \eqref{incompEuler2} as
\[
\Pi = {-}\Del^{-1} \Div \Div (\vc{v} \otimes \vc{v});
\]
whence it is uniformly bounded in $W^{1,q}(R^N)$ for any $1 < q < \infty$.

In accordance with \eqref{A1},
\[
\begin{split}
\int_0^\tau &\intOR{ \left[ \left< Y^\ep_{t,x}; \vr \vc{U} - \vc{m} \right> \cdot \left( \partial_t \Grad \Phi_\ep \right)
\right] }\dt = - \int_0^\tau \intOR{ \left< Y^\ep_{t,x}; \vc{m} \right> \partial_t \Grad \Phi_\ep } \dt\\
&+ \int_0^\tau \intOR{ \left< Y^\ep_{t,x}, \vr \right> \vc{v} \cdot \partial_t \Grad \Phi_\ep } \dt
+ \frac{1}{2} \int_0^\tau \intOR{ \left< Y^\ep_{t,x}, \vr \right> \cdot \partial_t |\Grad \Phi_\ep |^2  } \dt,
\end{split}
\]
where, as $\vc{v}$ is solenoidal,
\[
\begin{split}
\int_0^\tau &\intOR{ \left< Y^\ep_{t,x}, \vr \right> \vc{v} \cdot \partial_t \Grad \Phi_\ep } \dt \\ &= \ep
\int_0^\tau \intOR{ \left< Y^\ep_{t,x}, \frac{\vr - \Ov{\vr}}{\ep} \right> \vc{v} \cdot \partial_t \Grad \Phi_\ep } \dt =
- \frac{p'(\Ov{\vr})}{\Ov{\vr}} \int_0^\tau \intOR{ \left< Y^\ep_{t,x}, \frac{\vr - \Ov{\vr}}{\ep} \right> \vc{v} \cdot \Grad s_\ep } \dt.
\end{split}
\]
The rightmost expression tends to zero because of the dispersive estimates for $s_\ep$.

Similarly,
\[
\begin{split}
\frac{1}{2} \int_0^\tau &\intOR{ \left< Y^\ep_{t,x}, \vr \right> \cdot \partial_t |\Grad \Phi_\ep |^2  } \dt\\ &=
\ep \frac{1}{2} \int_0^\tau \intOR{ \left< Y^\ep_{t,x}, \frac{\vr - \Ov{\vr} }{\ep} \right> \cdot \partial_t |\Grad \Phi_\ep |^2  } \dt
+ \frac{\Ov{\vr}}{2} \left[ \intOR{ |\Grad \Phi_\ep |^2  } \right]^{t= \tau}_{t=0} \\ &=
-\frac{p'(\Ov{\vr})}{\Ov{\vr}} \int_0^\tau \intOR{ \left< Y^\ep_{t,x}, \frac{\vr - \Ov{\vr} }{\ep} \right> \cdot \Grad \Phi_\ep \cdot \Grad s_\ep   } \dt
+ \frac{\Ov{\vr}}{2} \left[ \intO{ |\Grad \Phi_\ep |^2  } \right]^{t= \tau}_{t=0},
\end{split}
\]
where the first term on the right-hand side vanishes for $\ep\to 0$ because of the dispersive estimates.

Consequently, the inequality (\ref{Il5}) can be recast in the form
\begin{equation} \label{Il6}
\begin{split}
&\mathcal{E}_\ep \left(\vr, \vc{m} \ \Big| \ \Ov{\vr} + \ep s_\ep , \vc{v} + \Grad \Phi_\ep \right)  (\tau)  + \mathcal{D}^\ep (\tau)
\leq \omega(\ep)
\\
&- \int_0^\tau \intOR{ \left< Y^\ep_{t,x}; \vc{m} \right> \partial_t \Grad \Phi_\ep } \dt + \frac{\Ov{\vr}}{2} \left[ \intOR{ |\Grad \Phi_\ep |^2  } \right]^{t= \tau}_{t=0}
\\
&
- \frac{1}{\ep^2} \int_0^\tau \intOR{ \left< Y^\ep_{t,x}; p(\vr) - p(\Ov{\vr}) \right> \Delta \Phi_\ep   } \dt \\
&+ \frac{1}{\ep} \int_0^\tau \intOR{ \left[ \left< Y^\ep_{t,x} ; r - \vr \right> P''(r) \partial_t s_\ep -
\left< Y^\ep_{t,x} ; {\vc{m}} \right> \cdot P''(r) \Grad s_\ep    \right] } \ \dt\\
& + c \int_0^\tau (1 +\xi(t)) \left[ \mathcal{E}
_\ep \left(\vr, \vc{m} \ \Big| \ \Ov{\vr} + \ep s_\ep , \vc{v} + \Grad \Phi_\ep \right) + \mathcal{D}^\ep \right] \dt
\end{split}
\end{equation}
where $\omega(\ep) \to 0$ as $\ep \to 0$.

\subsubsection{Step 3 - pressure estimates I}

We have
\[
\begin{split}
- \frac{1}{\ep} &\int_0^\tau \intOR{ \left< Y^\ep_{t,x}; \vc{m} \right> P''(r) \Grad s_\ep }\dt
\\ &= -  \int_0^\tau \intOR{ \left< Y^\ep_{t,x}; \vc{m} \right> \frac{P''(\Ov{\vr} + \ep s_\ep) - P''(\Ov{\vr}) }{\ep} \Grad s_\ep }\dt
- \frac{1}{\ep} \frac{p'(\Ov{\vr})}{\Ov{\vr}} \int_0^\tau \intOR{ \left< Y^\ep_{t,x}; \vc{m} \right> \cdot \Grad s_\ep }{\dt}\\
&=  -  \int_0^\tau \intOR{ \left< Y^\ep_{t,x}; \vc{m} \right> \frac{P''(\Ov{\vr} + \ep s_\ep) - P''(\Ov{\vr}) }{\ep} \Grad s_\ep }\dt
+ \int_0^\tau \intOR{ \left< Y^\ep_{t,x}; \vc{m} \right> \cdot \partial_t \Grad \Phi_\ep } \dt,
\end{split}
\]
where the first integral on the right-hand side vanished for $\ep \to 0$ by virtue of the dispersive estimates.

Next,
\[
\begin{split}
\frac{1}{\ep} \int_0^\tau &\intOR{ \left< Y^\ep_{t,x}; r - \vr  \right> P''(r) \partial_t s_\ep } \dt \\ &=
\frac{1}{2} \int_0^\tau \intOR{ P''(r) \partial_t |s_\ep|^2 } \dt +
\int_0^\tau \intOR{ \left< Y^\ep_{t,x}; \frac{ \Ov{\vr} - \vr }{\ep} \right> P''(r) \partial_t s_\ep } \dt \\&=
\left[ \frac{1}{2} \frac{p'(\Ov{\vr})}{\Ov{\vr}} \intOR{ |s_\ep |^2 } \right]_{t = 0}^{t = \tau} + \int_0^\tau \intOR{ \left< Y^\ep_{t,x}; \frac{ \Ov{\vr} - \vr }{\ep} \right> P''(r) \partial_t s_\ep } \dt\\&+
 \frac{\ep}{2} \int_0^\tau \intOR{ \frac{ P''(r) - P''(\Ov{\vr}) }{\ep}  \partial_t |s_\ep|^2 } \dt \\
&=
\left[ \frac{1}{2} \frac{p'(\Ov{\vr})}{\Ov{\vr}} \intOR{ |s_\ep |^2 } \right]_{t = 0}^{t = \tau} + \int_0^\tau \intOR{ \left< Y^\ep_{t,x}; \frac{ \Ov{\vr} - \vr }{\ep} \right> P''(r) \partial_t s_\ep } \dt\\& -
\Ov{\vr} \int_0^\tau \intOR{ \frac{ P''(r) - P''(\Ov{\vr}) }{\ep} s_\ep \Delta \Phi_\ep } \dt.
\end{split}
\]
Similarly to the above, the last integral is small in view of the dispersive estimates.

Summing up the previous observations with (\ref{Il6}) and using the acoustic energy balance (\ref{Il2bis}), we deduce from (\ref{Il6}) that
\begin{equation} \label{Il7}
\begin{split}
&\mathcal{E}_\ep \left(\vr, \vc{m} \ \Big| \ \Ov{\vr} + \ep s_\ep , \vc{v} + \Grad \Phi_\ep \right)  (\tau)  + \mathcal{D}^\ep(\tau)
\leq \omega(\ep)
\\
&
- \frac{1}{\ep^2} \int_0^\tau \intOR{ \left< Y^\ep_{t,x}; p(\vr) - p(\Ov{\vr}) \right> \Delta \Phi_\ep   } \dt +
\int_0^\tau \intOR{ \left< Y^\ep_{t,x}; \frac{\Ov{\vr} - \vr}{\ep} \right> P''(r) \partial_t s_\ep } \dt
\\
& + c\int_0^\tau (1+\xi(t)) \left[ \mathcal{E}_\ep \left(\vr, \vc{m} \ \Big| \ \Ov{\vr} + \ep s_\ep , \vc{v} + \Grad \Phi_\ep \right) + \mathcal{D}^\ep\right] \dt \end{split}
\end{equation}

\subsubsection{Step 4 - pressure estimates II}

We have
\[
\begin{split}
\int_0^\tau & \intOR{ \left< Y^\ep_{t,x}; \frac{\Ov{\vr} - \vr}{\ep} \right> P''(r) \partial_t s_\ep } \dt \\
&= \frac{p'(\Ov{\vr})}{\Ov{\vr}} \int_0^\tau \intOR{ \left< Y^\ep_{t,x}; \frac{\Ov{\vr} - \vr}{\ep} \right> \partial_t s_\ep } \dt
+ \ep \int_0^\tau  \intOR{ \left< Y^\ep_{t,x}; \frac{\Ov{\vr} - \vr}{\ep} \right> \frac{P''(r) - P''(\Ov{\vr}) }{\ep} \partial_t s_\ep } \dt
\\&= {p'(\Ov{\vr})} \int_0^\tau \intOR{ \left< Y^\ep_{t,x}; \frac{\vr - \Ov{\vr}}{\ep^2} \right> \Delta \Phi_\ep  } \dt
- \Ov{\vr} \int_0^\tau \intOR{\left< Y^\ep_{t,x}; \frac{\Ov{\vr} - \vr}{\ep} \right> \frac{P''(r) - P''(\Ov{\vr}) }{\ep} \Delta \Phi_\ep } \dt,
\end{split}
\]
where the last integral vanishes for $\ep \to 0$. Thus we obtain that
\[
\begin{split}
&\mathcal{E}_\ep \left(\vr, \vc{m} \ \Big| \ \Ov{\vr} + \ep s_\ep , \vc{v} + \Grad \Phi_\ep \right)  (\tau)  + \mathcal{D}^\ep (\tau) \leq \omega(\ep)
\\
&
- \frac{1}{\ep^2} \int_0^\tau \intOR{ \left< Y^\ep_{t,x}; \left[ p(\vr) - p'(\Ov{\vr})(\vr - \Ov{\vr}) - p(\Ov{\vr}) \right] \right> \Del \Phi_\ep   } \dt
\\
& + c \int_0^\tau (1 + \xi(t)) \left[ \mathcal{E}
_\ep \left(\vr, \vc{m} \ \Big| \ \Ov{\vr} + \ep s_\ep , \vc{v} + \Grad \Phi_\ep \right) + \mathcal{D}^\ep\right] \dt
\end{split}
\]
Consequently, using again the dispersive estimates (\ref{Il3}), together with the energy bounds,
we obtain the desired conclusion
\begin{equation} \label{IL8}
\begin{split}
&\mathcal{E}_\ep \left(\vr, \vc{m} \ \Big| \ \Ov{\vr} + \ep s_\ep , \vc{v} + \Grad \Phi_\ep \right)  (\tau)  + \mathcal{D}^\ep(\tau) \leq \omega(\ep)
\\
& + c \int_0^\tau (1 + \xi(t)) \left[ \mathcal{E}
_\ep \left(\vr, \vc{m} \ \Big| \ \Ov{\vr} + \ep s_\ep , \vc{v} + \Grad \Phi_\ep \right) + \mathcal{D}^\ep \right] \dt
\end{split}
\end{equation}
where $\xi \in L^1(0,T)$, and $\omega(\ep) \to 0$ as $\ep \to 0$.
Thus a direct application of Gronwall's lemma completes the proof of Theorem \ref{MRillp}.


\def\cprime{$'$} \def\ocirc#1{\ifmmode\setbox0=\hbox{$#1$}\dimen0=\ht0
  \advance\dimen0 by1pt\rlap{\hbox to\wd0{\hss\raise\dimen0
  \hbox{\hskip.2em$\scriptscriptstyle\circ$}\hss}}#1\else {\accent"17 #1}\fi}

\end{document}